\documentclass[12pt]{article}
\usepackage{fullpage}
\usepackage{epic}
\usepackage{eepic}
\usepackage{amsmath,amssymb,mathrsfs}
\usepackage{dsfont}
\usepackage[all]{xy}
\usepackage{enumerate}
\usepackage[T1]{fontenc}
\usepackage[latin1]{inputenc}
\usepackage{bbold}
\usepackage{soul}
\newtheorem{rem}{Remark}

\newtheorem{theo}{Theorem}
\newtheorem{lem}{Lemma}

\title{A short proof of the existence of supercuspidal representations for all reductive $p$-adic groups}
\author{Rapha\"{e}l Beuzart-Plessis\protect\footnote{This work was supported by the Gould Fund and by the National Science Foundation under agreement No. DMS-1128155. Any opinions, findings and conclusions or recommendations expressed in this material are those of the author and do not necessarily reflect the views of the National Science Foundation.}}
\begin{document}

\maketitle

\begin{abstract}
Let $G$ be a reductive $p$-adic group. We give a short proof of the fact that $G$ always admits supercuspidal complex representations. This result has already been established by A. Kret using the Deligne-Lusztig theory of representations of finite groups of Lie type. Our argument is of a different nature and is self-contained. It is based on the Harish-Chandra theory of cusp forms and it ultimately relies on the existence of elliptic maximal tori in $G$.
\end{abstract}

\vspace{7mm}

Let $p$ be a prime number and let $F$ be a $p$-adic field (i.e. a finite extension of $\mathds{Q}_p$). We denote by $\mathcal{O}$ the ring of integers of $F$ and we fix a uniformizer $\varpi\in \mathcal{O}$. We also denote by $val:F^\times\to \mathds{Z}$ the normalized valuation (i.e. $val(\varpi)=1$). Let $G$ be a connected reductive group defined over $F$. We will denote by $\mathfrak{g}$ the Lie algebra of $G$. A sentence like "Let $P=MN$ be a parabolic subgroup of $G$" will mean as usual that $P$ is a parabolic subgroup of $G$ defined over $F$, that $N$ is its unipotent radical and that $M$ is a Levi component of $P$ also defined over $F$. More generally, all subgroups of $G$ that we consider will be defined over $F$. We will also need to fix Haar measures on the various groups that we consider. The precise normalization of these Haar measures won't be important (unless we specify that they need to satisfy an explicit compatibility condition) and we will only make use of Haar measures on unimodular groups (e.g. $F$ points of reductive or unipotent groups) so that the distinction between left and right Haar measures is irrelevant here and will be dropped from the notations.

\vspace{2mm}

\begin{rem}\label{remark 1}
We exclude fields of positive characteristic because we will use in a crucial way the exponential map. If $G=GL_n$, we can use the map $X\mapsto Id+X$ instead and work over any non-archimedean local field. For classical groups, we could probably also replace the exponential map by some Cayley map and considerably weaken the characteristic assumption.
\end{rem} 

\vspace{2mm}

\noindent Recall that a smooth representation of $G(F)$ is a pair $(\pi,V_\pi)$ where $V_\pi$ is a complex vector space (usually infinite dimensional) and $\pi$ is a morphism $G(F)\to GL(V_\pi)$ such that for all vectors $v\in V_\pi$ the stabilizer $Stab_{G(F)}(v)$ of $v$ in $G(F)$ is an open subgroup. Let $(\pi,V_\pi)$ be a smooth representation of $G(F)$ and let $P=MN$ be a parabolic subgroup of $G$. The Jacquet module of $(\pi,V_\pi)$ with respect to $P$ is the space of coinvariants

$$V_{\pi,N}=V_\pi/V_\pi(N)$$

\noindent where $V_\pi(N)$ is the subspace of $V_\pi$ generated by the elements $v-\pi(n)v$ for all $v\in V_\pi$ and all $n\in N(F)$. It is the biggest quotient of $V_\pi$ on which $N(F)$ acts trivially. There is a natural linear action $\pi_N$ of $M(F)$ on $V_{\pi,N}$ and $(\pi_N,V_{\pi,N})$ is a smooth representation of $M(F)$. The functor $V_\pi\mapsto V_{\pi,N}$ is an exact functor from the category of smooth representations of $G(F)$ to the category of smooth representations of $M(F)$. Indeed, this follows rather easily from the following fact (cf \cite{Re} {\it Proposition III.2.9})

\vspace{5mm}

\hspace{5mm}(1) Let $\left(N(F)_k\right)_{k\geqslant 0}$ be an increasing sequence of compact-open subgroups of $N(F)$

\vspace{0.4mm}

\hspace{12mm} such that $N(F)=\bigcup_{k\geqslant 0} N(F)_k$ (such a sequence always exists). Then a vector

\vspace{0.4mm}

\hspace{12mm} $v\in V_\pi$ belongs to  $V_\pi(N)$ if and only if there exists $k\geqslant 0$ such that

$$\displaystyle \int_{N(F)_k}\pi(n)vdn=0$$

\vspace{5mm}

\noindent Let $(\pi,V_\pi)$ be an irreducible smooth representation of $G(F)$ (irreducible means that $V_\pi$ is non-trivial and that it has no nonzero proper $G(F)$-invariant subspace). The representation $(\pi,V_\pi)$ is said to be supercuspidal if for all proper parabolic subgroup $P=MN$ of $G$, the Jacquet module $V_{\pi,N}$ is zero. An equivalent condition is that the coefficients of $(\pi,V_\pi)$ are compactly supported modulo the center (cf \cite{Re} {\it Theorem VI.2.1}).

\vspace{2mm}

\noindent The purpose of this short note is to prove the following result.

\vspace{2mm}

\begin{theo}\label{theorem}
$G(F)$ admits irreducible supercuspidal representations.
\end{theo}

\begin{rem}
This theorem has already been proved by A.Kret (\cite{K}). The proof of Kret has the advantage of working in any characteristic (cf Remark \ref{remark 1}) and of being explicit (i.e. it exhibits a way to construct such supercuspidal representations by compact induction).  The proof of {\it loc. cit.} eventually relies on the Deligne-Lusztig theory of representations of finite groups of Lie type and so can hardly be considered as elementary. Although less complete and explicit than the results of {\it loc. cit.}, the proof presented here has the advantage of being short and (almost) self-contained using only elementary harmonic analysis arguments.
\end{rem}

\vspace{2mm}

\noindent We will deduce {\it Theorem \ref{theorem}} from the existence of nonzero compactly supported cusp forms, in the sense of Harish-Chandra, for the group $G(F)$. Before stating this existence result, we need to introduce some more definitions and notation. We will denote, as usual, by $C_c^\infty(G(F))$ the space of complex-valued functions on $G(F)$ that are smooth, i.e. locally constant, and compactly supported. We say that a function $f\in C_c^\infty(G(F))$ is a cusp form if for all proper parabolic subgroups $P=MN$ of $G$ we have

$$\displaystyle \int_{N(F)} f(xn)dn=0,\;\;\; \mbox{for all } x\in G(F)$$

\noindent (these functions are called supercusp forms in \cite{H-C}). We shall denote by $C_{c,cusp}^\infty(G(F))\subseteq C_c^\infty(G(F))$ the subspace of cusp forms. As we said, {\it Theorem \ref{theorem}} will follow from the following existence theorem.

\vspace{2mm}

\begin{theo}\label{proposition}
We have $C_{c,cusp}^\infty(G(F))\neq 0$.
\end{theo}

\vspace{2mm}

\noindent Proof that {\it Theorem \ref{proposition}} implies {\it Theorem \ref{theorem}}: Let us denote by $\rho$ the action of $G(F)$ on $C_c^\infty(G(F))$ given by right translation. Then, $\left(\rho, C_c^\infty(G(F))\right)$ is a smooth representation of $G(F)$. Moreover, it is easy to see that the subspace $C_{c,cusp}^\infty(G(F))\subseteq C_c^\infty(G(F))$ is $G(F)$-invariant. We claim the following:

\vspace{5mm}

\hspace{5mm} (2) For all proper parabolic subgroups $P=MN$ of $G$, we have

$$C_{c,cusp}^\infty(G(F))_N=0$$

\vspace{5mm}

\noindent Let $P=MN$ be a proper parabolic subgroup of $G$ and let us fix an increasing sequence $\left(N(F)_k\right)_{k\geqslant 0}$ of compact-open subgroups of $N(F)$ such that $N(F)=\bigcup_{k\geqslant 0} N(F)_k$. Let $f\in C_{c,cusp}^\infty(G(F))$. By (1), it suffices to show the existence of an integer $k\geqslant 0$ such that

$$\displaystyle \int_{N(F)_k} \rho(n)f dn=0$$

\noindent or what amounts to the same

$$(3)\;\;\; \displaystyle \int_{N(F)_k} f(xn)dn=0,\;\;\; \mbox{for all } x\in G(F)$$

\noindent Since $Supp(f)$ (the support of the function $f$) is compact, there exists $k\geqslant 0$ such that

$$(4)\;\;\;\displaystyle Supp(f)\cap \left[Supp(f)\left(N(F)\backslash N(F)_k\right)\right]=\emptyset$$

\noindent We now show that (3) is satisfied for such $k$. Let $x\in G(F)$. If $x\notin Supp(f)N(F)_k$, the term inside the integral (3) vanish identically and there is nothing to prove. Assume from now on that $x\in Supp(f)N(F)_k$. Up to translating $x$ by an element of $N(F)_k$, we may as well assume that $x\in Supp(f)$. Then, by (4) we have $xn\notin Supp(f)$ for all $n\in N(F)\backslash N(F)_k$. It follows that

$$\displaystyle \int_{N(F)_k} f(xn)dn=\int_{N(F)}f(xn)dn$$

\noindent But by definition of $C_{c,cusp}^\infty(G(F))$, this last integral vanishes. This proves (3) and ends the proof of (2).

\vspace{2mm}

\noindent We now show how to deduce from (2) that {\it Theorem \ref{proposition}} implies {\it Theorem \ref{theorem}}. Assume that {\it Theorem \ref{proposition}} is satisfied. Then, we can find $f\in C_{c,cusp}^\infty(G(F))$ which is nonzero. Denote by $V_f$ the $G(F)$-invariant subspace of $C_{c,cusp}^\infty(G(F))$ generated by $f$ and let $V\subseteq V_f$ be a maximal $G(F)$-invariant subspace among those not containing $f$ (Zorn's lemma). Then, $V_f/V$ is a smooth irreducible representation of $G(F)$. We claim that it is supercuspidal. Indeed, let $P=MN$ be a proper parabolic subgroup of $G$. By (2) and since the Jacquet module's functor is exact, we have $\left(V_f/V\right)_N=0$. Thus, $V_f/V$ is indeed a supercuspidal representation of $G(F)$ and this proves {\it Theorem \ref{theorem}}. $\blacksquare$

\vspace{3mm}

\noindent We are now left with proving {\it Theorem \ref{proposition}}. The strategy is to prove first an analog result on the Lie algebra and then lift it to the group by means of the exponential map. Let $C_c^\infty(\mathfrak{g}(F))$ be the space of complex-valued smooth and compactly supported functions on $\mathfrak{g}(F)$. We say that a function $\varphi\in C_c^\infty(\mathfrak{g}(F))$ is a cusp form if for all proper parabolic subgroup $P=MN$ of $G$ we have

$$\displaystyle \int_{\mathfrak{n}(F)}\varphi(X+N)dN=0,\;\;\; \mbox{for all } X\in \mathfrak{g}(F)$$

\noindent where $\mathfrak{n}(F)$ denotes the $F$-points of the Lie algebra of $N$. We will denote by $C_{c,cusp}^\infty(\mathfrak{g}(F))\subseteq C_c^\infty(\mathfrak{g}(F))$ the subspace of cusp forms. The analog of {\it Theorem \ref{proposition}} for the Lie algebra is the following lemma.

\vspace{2mm}

\begin{lem}\label{lemma}
We have $C_{c,cusp}^\infty(\mathfrak{g}(F))\neq 0$.
\end{lem}

\vspace{2mm}

\noindent Before proving this lemma, we first explain how it implies {\it Theorem \ref{proposition}}.

\vspace{2mm}

\noindent Proof that {\it Lemma \ref{lemma}} implies {\it Theorem \ref{proposition}}: Assume that {\it Lemma \ref{lemma}} holds. Then, we can find a nonzero function $\varphi\in C_{c,cusp}^\infty(\mathfrak{g}(F))$. The idea is to lift $\varphi$ to a function on $G(F)$ using the exponential map. Of course, the exponential map is not necessarily defined on the support of $\varphi$. Hence, we need first to scale the function $\varphi$ so that its support becomes small. Let us fix an element $\lambda\in F^\times$ that we will eventually assume to be sufficiently small. We define the function $\varphi_\lambda$ by

$$\varphi_\lambda(X)=\varphi(\lambda^{-1}X),\;\;\; X\in \mathfrak{g}(F)$$

\noindent We easily check that $\varphi_\lambda$ is still a cusp form. Recall that there exists an open neighborhood $\omega\subseteq \mathfrak{g}(F)$ of $0$ on which the exponential map $\exp$ is defined and such that it realizes an $F$-analytic isomorphism

$$\exp:\omega\simeq \Omega$$

\noindent where $\Omega=\exp(\omega)$. Since $Supp(\varphi_\lambda)=\lambda Supp(\varphi)$, for $\lambda$ sufficiently small we have 

$$Supp(\varphi_\lambda)\subseteq \omega$$

\noindent We henceforth assume that $\lambda$ is that sufficiently small. This allows us to define a function $f_\lambda$ on $G(F)$ by setting

$$f_\lambda(g)=\left\{
    \begin{array}{ll}
        \varphi_\lambda(X) & \mbox{ if } g=\exp(X) \mbox{ for some } X\in \omega \\
        0 & \mbox{ otherwise}
    \end{array}
\right.
$$

\noindent for all $g\in G(F)$. Note that we have $f_\lambda\in C_c^\infty(G(F))$ and obviously the function $f_\lambda$ is nonzero. Hence, we will be done if we can prove the following

\vspace{5mm}

\hspace{5mm} (5) If $\lambda$ is sufficiently small, the function $f_\lambda$ is a cusp form.

\vspace{5mm}

\noindent Let us denote by $\log:\Omega\to\omega$ the inverse of the exponential map. Then, by the Campbell-Hausdorff formula, it is easy to see that we can find an $\mathcal{O}$-lattice $L$ in the $F$-vector space $\mathfrak{g}(F)$ which is contained in $\omega$ and satisfies the following condition

$$(6)\;\;\; \displaystyle \log(e^Xe^Y)\in X+Y+\varpi^{val_L(X)+val_L(Y)}L$$

\noindent for all $X,Y\in L$, where we have set $val_L(X)=\sup\{k\in \mathds{Z}; \; X\in \varpi^kL\}$ for all $X\in\mathfrak{g}(F)$. For all integers $n\geqslant 0$, set $K_n=\exp(\varpi^n L)$. It is easy to infer from (6) that $K_n$ is an open-compact subgroup of $G(F)$ for all $n\geqslant 0$. Since $\varphi$ is smooth and compactly supported, there exists $n_0\geqslant 0$ such that translation by $\varpi^{n_0}L$ leaves $\varphi$ invariant. Also, since $\varphi$ is compactly supported, there exists $n_1\geqslant 0$ such that $Supp(\varphi)\subseteq \varpi^{-n_1}L$. We will show that (5) holds provided $val(\lambda)\geqslant 2n_1+n_0$. Assume this is so and set $n=val(\lambda)-n_1$. Then, we have

$$(7)\;\;\; Supp(\varphi_\lambda)=\lambda Supp(\varphi)\subseteq \lambda \varpi^{-n_1}L=\varpi^nL$$

\noindent Hence, it follows that

$$(8)\;\;\; Supp(f_\lambda)\subseteq K_n$$

\noindent Let $P=MN$ be a proper parabolic subgroup of $G$ and let $x\in G(F)$. Consider the integral

$$(9)\;\;\;\displaystyle \int_{N(F)} f_\lambda(xn)dn$$

\noindent If $xN(F)\cap K_n=\emptyset$, then by (8) the term inside the integral above vanishes identically and it follows that the integral is equal to zero. Assume from now on that $xK_n\cap N(F)\neq \emptyset$. Up to translating $x$ by an element of $N(F)$, we may assume that $x\in K_n$. Then, we may write $x=e^X$ for some $X\in \varpi^n L$. Using again (8), and since $K_n$ is a subgroup of $G(F)$, we see that the integral (9) is supported on $K_n\cap N(F)$. Thus, we have equalities

$$(10)\;\;\; \displaystyle \int_{N(F)} f_\lambda(xn)dn=\int_{K_n\cap N(F)} f_\lambda(e^Xn)dn$$

\noindent Set $L_N=L\cap \mathfrak{n}(F)$. Then, if we normalize measures correctly, the exponential map induces a measure preserving isomorphism $\varpi^n L_N\simeq K_n\cap N(F)$ so that we have

$$(11)\;\;\; \displaystyle \int_{K_n\cap N(F)} f_\lambda(e^Xn)dn=\int_{\varpi^n L_N} f_\lambda(e^Xe^N)dN=\int_{\varpi^n L} \varphi_\lambda(\log(e^Xe^N))dN$$

\noindent By (6), for all $N\in \varpi^n L_N$ we have

$$(12)\;\;\; \log(e^Xe^N)\in X+N+\varpi^{2n}L$$

\noindent Moreover, as $\varphi$ is invariant by translation by $\varpi^{n_0}L$, the function $\varphi_\lambda$ is invariant by translation by $\lambda \varpi^{n_0}L=\varpi^{n+n_1+n_0}L$ (recall that $n=val(\lambda)-n_1$). Since $val(\lambda)\geqslant 2n_1+n_0$, we also have $n\geqslant n_1+n_0$. Hence, the function $\varphi_\lambda$ is invariant by translation by $\varpi^{2n}L$ and so using (12), we deduce

$$\varphi_\lambda(\log(e^Xe^N))=\varphi_\lambda(X+N)$$

\noindent for all $N\in \varpi^n L_N$. From (10) and (11), it follows that

$$(13)\;\;\; \displaystyle \int_{N(F)} f_\lambda(xn) dn=\int_{\varpi^n L} \varphi_\lambda(X+N)dN$$

\noindent By (7) and since $X\in \varpi^n L$, the function $N\in \mathfrak{n}(F)\mapsto \varphi_\lambda(X+N)$ is supported on $\varpi^n L_N$. Consequently, we have

$$\displaystyle \int_{\varpi^n L} \varphi_\lambda(X+N)dN=\int_{\mathfrak{n}(F)}\varphi_\lambda(X+N)dN$$

\noindent As $\varphi_\lambda$ is a cusp form, this last integral vanishes. Hence, using (13) we see that the integral (9) is also zero. Since it is true for all $x\in G(F)$ and all proper parabolic subgroup $P=MN$ of $G$, this shows that $f_\lambda$ is a cusp form. Hence, (5) is indeed satisfied as soon as $val(\lambda)\geqslant 2n_1+n_0$ and this ends the proof that {\it Lemma \ref{lemma}} implies {\it Theorem \ref{proposition}}. $\blacksquare$

\vspace{3mm}

\noindent It only remains to establish {\it Lemma \ref{lemma}}. Recall that a maximal torus $T$ in $G$ is said to be elliptic if $A_T=A_G$, where $A_T$ and $A_G$ denote the maximal split subtori in $T$ and the center of $G$ respectively. The proof of {\it Lemma \ref{lemma}} will ultimately rely on the following existence result (cf \cite{PR} {\it Theorem 6.21}):

\vspace{2mm}

\begin{theo}\label{theorembis}
$G$ admits an elliptic maximal torus.
\end{theo}

\vspace{2mm}

\noindent Proof of {\it Lemma \ref{lemma}}: Let us fix a symmetric non-degenerate bilinear form $B$ on $\mathfrak{g}(F)$ which is $G(F)$-invariant. Such a bilinear form is easy to construct. On $\mathfrak{g}_{der}(F)$, the derived subalgebra of $\mathfrak{g}(F)$, we have the Killing form $B_{Kil}$ which is symmetric $G(F)$-invariant and non-degenerate. Hence, we may take $B=B_{\mathfrak{z}}\oplus B_{Kil}$ where $B_{\mathfrak{z}}$ is any non-degenerate symmetric bilinear form on $\mathfrak{z}_G(F)$, the center of $\mathfrak{g}(F)$. Let us also fix a non-trivial continuous additive character $\psi:F\to\mathds{C}^\times$. Using those, we can define a Fourier transform on $C_c^\infty(\mathfrak{g}(F))$ by

$$\displaystyle \widehat{\varphi}(Y)=\int_{\mathfrak{g}(F)}\varphi(X)\psi(B(X,Y))dX,\;\;\; \varphi\in C_c^\infty(\mathfrak{g}(F)),\; Y\in \mathfrak{g}(F)$$

\noindent Of course, this Fourier transform also depends on the choice of a Haar measure on $\mathfrak{g}(F)$. More generally, if $V$ is a subspace of $\mathfrak{g}(F)$ and $V^\perp$ denotes the orthogonal of $V$ with respect to $B$, we can also define a Fourier transform $C_c^\infty(V)\to C_c^\infty(\mathfrak{g}(F)/V^\perp)$, $\varphi\mapsto \widehat{\varphi}$, by setting

$$\displaystyle \widehat{\varphi}(Y)=\int_{V} \varphi(Y) \psi(B(X,Y))dY,\;\;\; X\in \mathfrak{g}(F)/V^\perp$$

\noindent where again we need to choose a Haar measure on $V$. It is easy to check that for compatible choices of Haar measures, the following diagram commutes

$$\xymatrix{
C_c^\infty(\mathfrak{g}(F)) \ar[d]^{res_V} \ar@{->}[r]^{FT} & C_c^\infty(\mathfrak{g}(F)) \ar[d]^{\int_{V^\perp}}\\
C_c^\infty(V) \ar@{->}[r]^{FT} & C_c^\infty(\mathfrak{g}(F)/V^\perp)}$$

\noindent where the horizontal arrows are Fourier transforms, the left vertical arrow is given by restriction to $V$ and the right vertical arrow is given by integration over the cosets of $V^\perp$. For $P=MN$ a parabolic subgroup of $G$, we have $\mathfrak{p}(F)^\perp=\mathfrak{n}(F)$ (where $\mathfrak{p}$ stands for the Lie algebra of $P$). The commutation of the above diagram in this particular case gives us (for some compatible choices of Haar measures) the following formula

$$(14)\;\;\; \displaystyle \int_{\mathfrak{n}(F)} \widehat{\varphi}(X+N)dN=\int_{\mathfrak{p}(F)} \varphi(Y) \psi(B(X,Y))dY$$

\noindent for all $\varphi\in C_c^\infty(\mathfrak{g}(F))$ and all $X\in \mathfrak{g}(F)$.

\vspace{2mm}

\noindent Let $T_{ell}$ be an elliptic maximal torus of $G$ whose existence is guaranteed by {\it Theorem \ref{theorembis}}. Let $\mathfrak{t}_{ell}$ be its Lie algebra and set $\mathfrak{t}_{ell,reg}=\mathfrak{t}_{ell}\cap \mathfrak{g}_{reg}$ for the subset of $G$-regular elements in $\mathfrak{t}_{ell}$. Denote by $\mathfrak{t}_{ell,reg}(F)^G$ the subset of elements in $\mathfrak{g}_{reg}(F)$ that are $G(F)$-conjugated to an element of $\mathfrak{t}_{ell,reg}(F)$. Then, $\mathfrak{t}_{ell,reg}(F)^G$ is an open subset of $\mathfrak{g}(F)$ (since the map $T_{ell}(F)\backslash G(F)\times \mathfrak{t}_{ell,reg}(F)\to \mathfrak{g}(F)$, $(g,X)\mapsto g^{-1}Xg$, is everywhere submersive). In particular, we can certainly find a non-zero function $\varphi\in C_c^\infty(\mathfrak{g}(F))$ whose support is contained in $\mathfrak{t}_{ell,reg}(F)^G$. Let us fix such a function $\varphi$. We claim the following:

\vspace{5mm}

\begin{center}
(15) The function $\widehat{\varphi}$ is a cusp form.
\end{center}

\hspace{5mm}

\noindent Indeed, let $P=MN$ be a proper parabolic subgroup of $G$ and let $X\in \mathfrak{g}(F)$. Then, we need to see that the following integral

$$\displaystyle \int_{\mathfrak{n}(F)} \widehat{\varphi}(X+N)dN$$

\noindent is zero. By (14), this integral is equal to

$$\displaystyle \int_{\mathfrak{p}(F)} \varphi(Y)\psi(B(X,Y))dY$$

\noindent Hence, we only need to show that $Supp(\varphi)\cap \mathfrak{p}(F)=\emptyset$. By definition of $\varphi$, it even suffices to see that $\mathfrak{t}_{ell,reg}(F)^G\cap \mathfrak{p}(F)=\emptyset$. But this follows immediately from the fact that $P$ being proper, it does not contain any elliptic maximal torus of $G$. $\blacksquare$

\bigskip

\hspace{3mm} Institute for Advanced Study, Princeton, NJ USA

\hspace{3mm} email address: rbeuzart@gmail.com

\end{document}